\documentclass[12pt, 14paper,reqno]{amsart}
\usepackage{amsmath,amsfonts,amssymb}
\usepackage[breaklinks]{hyperref}
\usepackage{graphicx}
\usepackage{longtable}

\theoremstyle{plain}
\newtheorem{thm}{Theorem}[section]
\newtheorem{lem}{Lemma}[section]
\newtheorem{cor}{Corollary}[section]
\newtheorem{prop}{Proposition}[section]

\theoremstyle{proof}
\makeatletter
\@namedef{subjclassname@2020}{%
  \textup{2020} Mathematics Subject Classification}
\makeatother

\numberwithin{equation}{section}

\begin{document} 
\title[On an exponential Diophantine equation]{Complete solutions of a Lebesgue-Ramanujan-Nagell type equation}
\author{Priyanka Baruah, Anup Das and Azizul Hoque}
\address{Department of Mathematics, Gauhati University, Guwahati-781014, Assam, India.}
\email{baruah067@gmail.com}
\address{Department of Mathematics, Gauhati University, Guwahati-781014, Assam, India.}
\email{anup@gauhati.ac.in}
\address{Department of Mathematics, Faculty of Science, Rangapara College, Rangapara, Sonitpur-784505, Assam, India.}
\email{ahoque.ms@gmail.com}
\keywords{Diophantine equation, Lehmer sequence, Elliptic curve, Quartic curve, S-Integral Points}
\subjclass[2020] {Primary: 11D61, 11D41, Secondary: 11Y50}
\date{\today}
\maketitle

\begin{abstract}
We consider the Lebesgue-Ramanujan-Nagell type equation $x^2+5^a13^b17^c=2^m y^n$, where $a,b,c, m\geq 0, n \geq 3$ and $x, y\geq 1$ are unknown integers with $\gcd(x,y)=1$. We determine all integer solutions to the above equation. The proof depends on the classical results of Bilu, Hanrot and Voutier on primitive divisors in Lehmer sequences, and finding all $S$-integral points on a class of elliptic curves.  
\end{abstract}
\section{Introduction}
The Lebesgue-Ramanujan-Nagell type equation
\begin{equation}\label{eqi1}
x^2 + D^m =\lambda y^n, ~\lambda=1,2,4,
\end{equation} in integer unknowns $x, y, m\geq 1$ and $n \geq 3$, has a long and distinguished history. The first result concerning the solutions of \eqref{eqi1} was due to Lebesgue \cite{LE50}, who proved that \eqref{eqi1} has no solutions when $\lambda=D=1$ and $y>1$.  Later, many authors become interested in this equation and thus there are good amount of research concerning the solutions of \eqref{eqi1}. We direct the readers to  the beautiful survey \cite{LS20} for further information. Several authors studied some generalizations of \eqref{eqi1} in \cite{AA02, BU01, CH22, CHS21, DGS21, HO20}.  

Recently, many authors become interested to find the integer solutions of the Lebesgue-Ramanujan-Nagell type equation
$$x^2+p_1^{a_1}p_2^{a_2}\cdots p_k^{a_k}=y^n,~~x,y\geq 1, \gcd(x,y)=1, a_1,a_2,\cdots, a_k\geq 0, n\geq 3,$$
where $p_1, p_2, \cdots, p_k$ are distinct primes with $k\geq 2$. There are many results concerning the integer solutions of this equation, but we refer to the very recent papers \cite{AL09, AZ20, BHS19, CHS20, CHS-21, DA11, DE17,  ZL15}. From existing results, it is quite natural to consider Diophantine equations similar to the above one, where the right side is a product of an unknown integer with unknown exponent and a known prime with unknown exponent. 

Here, we consider the Diophantine equation
\begin{equation}\label{eqn1}
x^2+5^a13^b17^c=2^m y^n,~~ x\geq 1, y>1, \gcd(x, y)=1, a, b,c, m\geq 0, n\geq 3,
\end{equation}
and we find all its the integer solutions $(x, y, a,b,c,m, n)$. It is noted that by reading \eqref{eqn1} modulo $4$, we see that $x^2+1\equiv 0\pmod 4$ which does not satisfy for any integer $x$. Therefore \eqref{eqn1} has no solution when $m\geq 2$, and  thus we consider \eqref{eqn1} for $m=0, 1$. More precisely, we prove the following:

\begin{thm}\label{thm}
If $n\ne 3,4,6,12$, then \eqref{eqn1} has no integer solutions. In case of $n=3,4,6,12$, the integer solutions $(x,y,a,b,c,m)$ are given below. 
\begin{itemize}  
\item[(i)] For $n=3$,  $(x,y,a,b,c,m)$ are given in  Table \ref{Tn}.
\item[(ii)] For $n=4$, $(x,y,a,b,c,m)$ are given in  Table \ref{Tm4}.
\item[(iii)] For $n=6$, $(x,y,a,b,c,m)=(716,9,1,1,2,0)$.
\item[(iv)] For $n=12$,  $(x,y,a,b,c,m)=(716,3,1,1,2,0)$.
\end{itemize}
\end{thm}

\subsection*{Remarks} We mention some earlier results which can be retrieved from Theorem \ref{thm}.
\begin{itemize}
\item[(i)] For $m\geq 2$, reducing \eqref{eqn1} modulo $4$, one can see that it has no integer solutions.
\item[(ii)] Abu Muriefah and Arif proved that the Diophantine equation $x^2+5^a=y^n, n\geq 3, x\geq 1, y>1, \gcd(x,y)=1$, has no integer solutions when $a$ is odd (see,  \cite[Theorem]{AA99}). Later, Tao completely solved it in \cite{TA09} and proved that it has no integer solutions. We can get these results from our theorem. In fact, our theorem shows that $(x,y,a,m, n)=(239,13,0,1,4), (7,3,1,1,3), (99,17,2,1,3)$ are the only integer solutions of the Diophantine equation $$x^2+5^a=2^my^n, n\geq 3, x\geq 1, y>1, \gcd(x,y)=1.$$  

\item[(iii)] It follows from \cite[Theorem 1.1]{LT09} that $(x, y, b, n) = (70, 17, 1, 3)$ is the only integer solution of the Diophantine equation $x^2+13^b=y^n~~ (b, x,y\geq 1, \gcd(x,y)=1, n\geq 3)$. A consequence of Theorem \ref{thm} extends this result to the Diophantine equation $x^2+13^b=2^my^n~~ (x,y\geq 1, \gcd(x,y)=1, b, m\geq 0, n\geq 3)$. In this case, the only integer solutions are $(x,y,b,m, n)=(70,17,1,0,3), (9,5,2,1,3), (239,13,0,1,4)$.

\item[(iv)] In \cite{AL08}, Abu Muriefah et al. proved that the Diophantine equation $x^2+5^a13^b=y^n, n\geq 3, x\geq 1, y>1, \gcd(x,y)=1$, has no integer solutions, except $(x,y,a,b) = (70,17,0,1), (142,29,2,2), (4,3,1,1)$. We can get this result from our theorem. Our theorem also confirms that the integer solution of the Diophantine equation $x^2+5^a13^b=2y^n ~~(n\geq 3, x\geq 1, y>1, \gcd(x,y)=1)$, are $(x,y,a,b,n)=(239,13,0,0,4), (9,5,0,2,3), (7,3,1,0,3), (99,17,2,0,3), 
 (19,7,2,1,3), (253,73,2,4,3), $ $(79137, 1463, 2,3,3),  (188000497, 260473, 8,4,3)$. 

\item[(v)] Pink and R\'abai completely solved the  Diophantine equation $x^2+5^a17^c=y^n~~ (x, y\geq 1, \gcd(x,y)=1, a,c\geq 0, n\geq 3)$ in \cite{PR11}. However, our theorem gives all the integer solutions of an extension of this equation, namely $x^2+5^a17^c=2^my^n~~ (x, y\geq 1, \gcd(x,y)=1, a,c\geq 0, n\geq 3)$.
\end{itemize}

Theorem \ref{thm} yields the following straightforward corollary. In case of $m=0$, this corollary follows from the work of Gou and Wang \cite{GW12}.
\begin{cor}
The Diophantine equation $$x^2+17^k=2^my^n,~~x\geq 1, y>1, \gcd(x,y)=1, k, m\geq 0, n\geq 3,$$ has no integer solutions, except $(x,y,k,m,n)=(8,3,1, 0,4), (31,5,2,1,4), (239,13,0,1,4)$.
\end{cor}

The next corollary immediately follows from Theorem \ref{thm}. 
\begin{cor}
The only integer solutions of the Diophantine equation $$x^2+13^k17^\ell=2^my^n,~~x\geq1, y>1, \gcd(x,y)=1, k, \ell, m\geq 0, n\geq 3$$ are $(x,y,k,\ell, m,n)=(70, 17,1,0,0,3), (9,5,2,0,1,3), (8, 3,0,1,0,4), (31,5,0,2,1,4),$\\ $ (239,13,0,0,1,4)$.
\end{cor}

We organize this article as follows. In \S\ref{S2}, we deal with the exponent $n$ satisfying  $4\mid n$. In this case, we transform \eqref{eqn1} into quartic curves, and thus the problem is reduced to finding all $\{5,13,17\}$-integral points on these curves.  Recall that for a finite set of prime numbers $S$, an $S$-integer is a rational number $r/s$ with coprime integers $r$ and $s>0$ such that any prime factor of $s$ lies in $S$. We treat \eqref{eqn1} in \S\ref{S3} for prime exponent $n\geq 3$. For $n\geq 5$ with $n\ne 7$, we apply the result of Bilu, Hanrot and Voutier \cite{BH01} concerning the existence of primitive divisors in Lehmer sequences. In case of $n=7$, we first use  some criteria for the existence of primitive divisors in Lehmer sequences to handle some cases of \eqref{eqn1}. For the remaining cases, we somehow transform them into elliptic curves. Analogously, we transform \eqref{eqn1} into elliptic curves for $n=3$. Then we solve the problem by finding all $\{5,13,17\}$-integral points on these elliptic curves. In \S\ref{S4}, we summarize the proof of Theorem \ref{thm}. All the computations are done using MAGMA \cite{BCP97}.

\section{The case: $4\mid n$}\label{S2}
Here, we prove the following proposition.
\begin{prop}\label{propm4}
If $n$ is a multiple of $4$, then all integer solutions of \eqref{eqn1} are given in Table \ref{Tm4}. 
\end{prop}\vspace{-6mm}
			
\begin{center}
{\small
\begin{longtable}{ccccccc| ccccccc}
\caption{All the solutions of \eqref{eqn1} when $4\mid n$} \label{Tm4} \\

\hline \multicolumn{1}{c}{$x$} & \multicolumn{1}{c}{$y$} & \multicolumn{1}{c}{$a$}& \multicolumn{1}{c}{$b$} & \multicolumn{1}{c}{$c$}& \multicolumn{1}{c}{$m$} &\multicolumn{1}{c}{$n$} &\multicolumn{1}{|c}{$x$} & \multicolumn{1}{c}{$y$} &\multicolumn{1}{c}{$a$} &\multicolumn{1}{c}{$b$} & \multicolumn{1}{c}{$c$}& \multicolumn{1}{c}{$m$} & \multicolumn{1}{c}{$n$}\\ \hline 
\endfirsthead

\multicolumn{14}{c}%
{{\bfseries \tablename\ \thetable{} -- continued from previous page}} \\
\hline \multicolumn{1}{c}{$x$} & \multicolumn{1}{c}{$y$} & \multicolumn{1}{c}{$a$}& \multicolumn{1}{c}{$b$} & \multicolumn{1}{c}{$c$}& \multicolumn{1}{c}{$m$} &\multicolumn{1}{c}{$n$} &\multicolumn{1}{|c}{$x$} & \multicolumn{1}{c}{$y$} &\multicolumn{1}{c}{$a$} &\multicolumn{1}{c}{$b$} & \multicolumn{1}{c}{$c$}& \multicolumn{1}{c}{$m$} & \multicolumn{1}{c}{$n$}\\ \hline 
\endhead
\hline \multicolumn{14}{c}{{Continued on next page}} \\ \hline
\endfoot
\hline 
\endlastfoot
8&3&0&0&1&0&4&
4&3&1&1&0&0&4\\
26556&163&5&1&1&0&4&
36&7&1&1&1&0&4\\
716&27&1&1&2&0&4&
716&3&1&1&2&0&12\\
239&13&0&0&0&1&4&
31&5&0&0&2&1&4
\end{longtable}}
\end{center}
\vspace{-14mm}

\begin{proof}
Assume that $n=4t$, where $t\geq 1$ is an integer. Then \eqref{eqn1} can be written as 
\begin{equation}\label{eqm41}
x^2+5^a13^b17^c=2^m\left(y^t\right)^4,~ x\geq 1, y>1, \gcd(x, y)=1, a,b,c, m\geq 0, t\geq 1.\
\end{equation}
\vspace{1mm}
Recall that \eqref{eqm41} has no integer solution when $m\geq 2$. 
Let $a\equiv a_1\pmod 4, b\equiv b_1\pmod 4$ and $ c\equiv c_1\pmod 4$. Then \eqref{eqm41} can be written as
$$x^2+5^{a_1}13^{b_1}17^{c_1}z^4=2^m\left(y^t\right)^4,$$
where $5^a13^b17^c=5^{a_1}13^{b_1}17^{c_1}z^4$. 
This can be written as
\begin{equation}\label{eqm42}
X^2=2^mY^4-5^{a_1}13^{b_1}17^{c_1},
\end{equation}
where $X=x/z^2$ and $Y=y^t/z$. Now the  problem of finding integer solutions of \eqref{eqm41} is transformed to finding all $\{5,13,17\}$-integer points on the $128$ quartic curves defined by \eqref{eqm42}. Here, we use MAGMA \cite{BCP97} subroutine \texttt{SIntegralLjunggrenPoints} to determine all $\{5,13,17\}$-integer points on these elliptic curves. Note that we avoid $\{5,13,17\}$-integer points with $XY = 0$ as they yield to $xy = 0$. Also taking into account that $\gcd(x, y) = 1$, we don't consider $\{5,13,17\}$-integer points such that the numerators of $X$ and $Y$ are not coprime. We finally get only $8$ integer points $(X, Y)$ with $XY\ne 0$ and the numerators of $X$ and $Y$ are coprime. We then use the relations
$$X=\frac{x}{z^2}, Y=\frac{y^t}{z} \text{  and }5^a13^b17^c=5^{a_1}13^{b_1}17^{c_1}z^4,$$
to find the integer solutions $(x,y,a,b,c,m, n)$, which are listed in Table \ref{Tm4}. 
\end{proof}
\section{The case: $n\geq 3$ is prime}\label{S3}
We rewrite \eqref{eqn1} by changing $n$ to $p$ to emphasize that the exponent is prime:
\begin{equation}\label{eqp1}
x^2+5^a13^b17^c=2^my^p,~ x,\geq 1, y>1, \gcd(x, y)=1, a,b,c\geq 0, p\geq 2, m=0,1.
\end{equation}
\begin{prop}\label{propp} The equation \eqref{eqp1} has no integer solutions for $p>3$. When $p=3$, its integer solutions are given by  
$(x,y,a,b,c,m)\in\mathfrak{S}$, where
\begin{align*}
\mathfrak{S}:=&\{(70,17,0,1,0,0), (716,81,1,1,2,0), (94,21,2,0,1,0),(142,29,2,2,0,0),
(2034,161,3,\\
&0,2,0), (9,5,0,2,0,1), (7,3,1,0,0,1), (99,17,2,0,0,1), 
(63,13,2,0,1,1), (19,7,2,1,0,1), \\
& (33,7,2,2,1,1), (118699, 1917,2,2,1,1), (79137, 1463, 2,3,0,1), (253,73,2,4,0,1), \\
&(188000497, 260473, 8,4,0,1), (267689, 3297, 2,2,3,1), (336049, 4317, 10,0,3,1), \\ 
&(17127, 553, 6,2,1,1)\}.
\end{align*}

\end{prop}
Before proceeding further, we need to recall some results and to fix some notations. The following lemma follows from \cite[Corollary 3.1]{YU05}; however, the idea goes back to the work of Ljunggren \cite{LJ66}.

\begin{lem}\label{lemYU}
Let $d~ (\ne 3)$ be square-free positive integer. If $n\geq 3$ is an odd integer coprime to  $h(-d)$, the class number of $\mathbb{Q}(\sqrt{-d})$, then all integer solutions $(X,Y,Z)$ of the Diophantine equation 
\begin{equation*}\label{eqYU}
X^2+dY^2=2^m Z^n,~~ X,Y\geq 1, \gcd(X, dY)=1, m=0,1,
\end{equation*} 
can be expressed as 
$$\frac{X+Y\sqrt{-d}}{\sqrt{2^m}}=\varepsilon_1\left(\frac{u+\varepsilon_2 v\sqrt{-d}}{\sqrt{2^m}}\right)^n,$$
where $u$ and $v$ are positive integers satisfying $2^mZ=u^2+dv^2$ and $\gcd(u, dv)=1$, and $\varepsilon_1,\varepsilon_2\in\{-1, 1\}$.
\end{lem}

Assume that $F_k$ (resp. $L_k$) denote the $k$-th Fibonacci (resp. Lucas) number defined by $F_0=0, F_1=1, F_k=F_{k-1}+F_{k-2}$ (resp. $L_0=2, L_1=1, L_k=L_{k-1}+L_{k-2}$) for $k\geq 2$. Then the following lemma follows from Theorems 1,2,3 and 4 in \cite{CO64}. 
\begin{lem}\label{lemCO}
Let $F_k$ and $L_k$ be the $k$-th Fibonacci and Lucas numbers, respectively. Then
\begin{itemize}
\item[(i)] if $F_k = x^2$, then $(k,x)=(0, 0), (1, \pm 1),(2, \pm 1), (12, \pm 12)$;
\item[(iii)] if $F_k =2 x^2$, then $(k,x)=(0, 0), (3, \pm 1), (6,\pm 2)$;
\item[(iii)] if $L_k = x^2$, then $(k,x)=(1, \pm 1), (3, \pm 2)$;
\item[(iv)] if $L_k = 2x^2$, then $(k,x)=(0, \pm 1), (6, \pm 3)$.
\end{itemize}
\end{lem}
A pair $(\alpha_1, \alpha_2)$ of algebraic integers is said to be a Lehmer pair if $(\alpha_1 + \alpha_2)^2$ and $\alpha_1\alpha_2$ are two non-zero coprime rational integers, and $\alpha_1/\alpha_2$ is not a root of unity. Also for a positive integer $n$, the Lehmer number corresponds to the pair $(\alpha_1, \alpha_2)$ is defined as 
$$\mathfrak{L}_n(\alpha_1, \alpha_2)=\begin{cases}
\dfrac{\alpha_1^n-\alpha_2^n}{\alpha_1-\alpha_2} & \text{ if } n \text{ is odd}, \vspace{2mm}\\
\dfrac{\alpha_1^n-\alpha_2^n}{\alpha_1^2-\alpha_2^2} & \text{ if } n \text{ is even}.
\end{cases}$$
Note that all Lehmer numbers are non-zero rational integers. A prime divisor $p$ of $\mathfrak{L}_n(\alpha_1,\alpha_2)$ is primitive if $p\nmid(\alpha_1^2-\alpha_2^2)^2
\mathfrak{L}_1(\alpha_1, \alpha_2) \mathfrak{L}_2(\alpha_2, \alpha_2) \cdots \mathfrak{L}_{n-1}(\alpha_1, \alpha_2)$. 
Further,
$\left( (\alpha_1+\alpha_2)^2, (\alpha_1-\alpha_2)^2\right)$ is known as the parameters of the Lehmer pair $(\alpha_1, \alpha_2)$.

\begin{proof}[{\bf{Proof of Proposition \ref{propp}}}]
We first rewrite \eqref{eqp1} as follows:
\begin{equation}\label{eqp2}
x^2+dz^2=2^my^p,~ x,y\geq 1, \gcd(x, y)=1, a,b,c\geq 0, p\geq 2, m=0,1,
\end{equation}
where $d\in\{1, 17, 13, 221, 5,85,65,1105\}$ and $z=5^{a_1}13^{b_1}17^{c_1}$ for some integers $a_1,b_1,c_1\geq 0$.
Using MAGMA, we see that $h(-d)\in\{1,2,4,8,16\}$, and thus $p\nmid h(-d)$. As $\gcd(x,y)=1$, so that $\gcd(x,dz)=1$. Therefore by Lemma \ref{lemYU}, we have (from \eqref{eqp2})
\begin{equation}\label{eqp3}
\frac{x+z\sqrt{-d}}{\sqrt{2^m}}=\varepsilon_1\left(\frac{u+\varepsilon_2 v\sqrt{-d}}{\sqrt{2^m}}\right)^p,
\end{equation}
where $u$ and $v$ are positive integers satisfying $2^my=u^2+dv^2$ and $\gcd(u, dv)=1$, and $\varepsilon_1,\varepsilon_2\in\{-1, 1\}$. Note that $2\nmid duvy$. 

We define, 
$$\alpha:=\dfrac{u+\varepsilon_2 v\sqrt{-d}}{\varepsilon_1\sqrt{2^m}}.$$
Then $\alpha$ and its conjugate $\bar{\alpha}$ are algebraic integers such that $\gcd((\alpha+\bar{\alpha})^2,\alpha\bar{\alpha})=1$. It is easy to see that $\alpha/\bar{\alpha}$ satisfies $$yZ^2-2^{1-m}(u^2-v^2d)T+y=0.$$ This shows that $\alpha/\bar{\alpha}$ is not a root of unity as $\gcd(2^{m-1}(u^2-v^2d), y)=1$. Thus $(\alpha, \bar{\alpha})$ is a Lehmer pair and $(2^{2-m}u^2, -2^{2-m}v^2d)$ is the corresponding parameter. 

Let $\mathfrak{L}_n$ be the Lehmer number corresponding to the Lehmer pair $(\alpha, \bar{\alpha})$. Then 
\begin{equation}\label{eqp4}
|\mathfrak{L}_p(\alpha, \bar{\alpha})|=\frac{z}{v}=\frac{5^{a_1}13^{b_1}17^{c_1}}{v}.
\end{equation}  
We divide the remaining part of the proof into several parts depending of the values of $p$.

\subsection*{Case I: When $p>7$} \label{SS1} 
Assume that $q$ is a primitive divisor of $\mathfrak{L}_p(\alpha, \bar{\alpha})$. Then by \eqref{eqp4}, $q\in\{5,13,17\}$. We now utilize the fact that any primitive divisor of $\mathfrak{L}_p(\alpha, \bar{\alpha})$ is congruent to $\pm 1$ modulo $p$, to conclude that none of these values of $q$ is primitive divisor of $\mathfrak{L}_p(\alpha, \bar{\alpha})$.  This contradicts to a consequence of the Primitive Divisor Theorem for Lehmer sequences which states that, if $p\geq 3$, then $\mathfrak{L}_p(\alpha, \bar{\alpha})$ has a primitive prime divisor except for finitely many pairs $(\alpha, \bar{\alpha})$. The Lehmer sequences correspond to these exceptional pairs $(\alpha, \bar{\alpha})$ are given in \cite[Tables 2 and 4]{BH01} in terms of their parameters $((\alpha+\bar{\alpha})^2,(\alpha-\bar{\alpha})^2)$.
Since $(2^{2-m}u^2, -2^{2-m}v^2d)$ are the parameters, so that $p=13$ and $(2^{2-m}u^2, 2^{2-m}v^2d)=(1,7)$, which is not possible as $m=0, 1$. This concludes that \eqref{eqp1} has no integer solutions when $p>7$. 

\subsection*{Case II: When $p=7$} As in Case I, if $q$ is a primitive divisor of $\mathfrak{L}_7(\alpha, \bar{\alpha})$, then only possibility is $q=13$. If $b_1=0$, then $\mathfrak{L}_7(\alpha, \bar{\alpha})$ has no primitive divisors, and hence as in Case I by \cite[Table 2]{BH01}, we have $(2^{2-m}u^2, 2^{2-m}v^2d)=(1,7), (1,19), (3,5),(5,7), (13,3),(14,22)$. These are not possible as $m=0,1$. Therefore \eqref{eqp1} has no solutions. Analogously, we can conclude that $\eqref{eqp1}$ has no solutions if $a_1=b_1=c_1=0$.

We now consider $b_1\geq 1$. Note that 
\begin{equation}\label{eqp5}
(\alpha^2-\bar{\alpha}^2)^2=-2^{4-2m}u^2v^2d.
\end{equation}
Thus if $13\mid vd$ then $13$ is not a primitive divisor, and hence as before \eqref{eqp1} has no solutions. Therefore $d=1,5,17,85$ and $13\nmid v$. Also using the fact that for a primitive divisor $q$, the sign of $q\equiv \pm 1\pmod p$ coincides with that of the Legendre symbol $\left(\frac{-4d}{q}\right)$, we get $d=5,85$. Equating imaginary parts in \eqref{eqp3}, we get
\begin{equation}\label{eqp6}
2^{3m}5^{a_1}13^{b_1}17^{c_1}=\varepsilon v(7u^6-35u^4v^2d+21u^2v^4d^2-v^6d^3),
\end{equation}
where $\varepsilon=\varepsilon_1\varepsilon_2=\pm 1$.
Since $\gcd(u,v)=1$ and $2,13\nmid v$, so that \eqref{eqp6} gives $$v=1, 5^{a_1},17^{c_1}, 5^{a_1}17^{c_1}.$$

For $v=1$, \eqref{eqp6} becomes 
$$2^{3m} 5^{a_1}13^{b_1}17^{c_1}=\varepsilon(7u^6-35u^4d+21u^2d^2-d^3).
$$
We can rewrite this equation as
\begin{equation}\label{eqp7}
DY^2=7X^3-35dX^2+21d^2X-d^3,
\end{equation}
where $X=u^2, Y=2^{m}5^{a_2}13^{b_2}17^{c_2}, a_2=\lfloor{a_1/2}\rfloor, b_2=\lfloor{b_1/2}\rfloor, c_2=\lfloor{c_1/2}\rfloor$ and $D=\varepsilon 2^m5^i13^j17^k$ with $i, j, k\in\{0,1\}$.
We now multiply both sides of \eqref{eqp7} by $7^2D^3$ and then rewrite as follows:
\begin{equation}\label{eqp8}
V^2=U^3-35dDU+147d^2D^2U-49d^3D^3,
\end{equation}
where $U=7DX$ and $V=7D^2Y$. Here, we use \texttt{IntegralPoints} subroutine of MAGMA \cite{BCP97} to compute all integral points on the $64$ elliptic curves defined by \eqref{eqp8}. We then apply the relations $U=7DX, V=7D^2Y, X=u^2, Y=2^{m}5^{a_2}13^{b_2}17^{c_2}, a_2=\lfloor{a_1/2}\rfloor, b_2=\lfloor{b_1/2}\rfloor, c_2=\lfloor{c_1/2}\rfloor$ and $D=\varepsilon 2^m5^i13^j17^k$ with $i, j, k\in\{0,1\}$. 
We check that none of these integral points leads to an integer solution of \eqref{eqp1}.

We now consider $v=5^{a_1}$ with $a_1\geq 1$. Then \eqref{eqp6} becomes 
$$2^{3m}13^{b_1}17^{c_1}=\varepsilon(7u^6-35u^4v^2d+21u^2v^4d^2-v^6d^3).
$$
Dividing both sides of this equation by $v^6$, we obtain the following elliptic curves
\begin{equation}\label{eqp9}
DY^2=7X^3-35dX^2+21d^2X-d^3,
\end{equation}
where $X=u^2/v^2, Y=2^{m}13^{b_2}17^{c_2}/v^3, b_2=\lfloor{b_1/2}\rfloor, c_2=\lfloor{c_1/2}\rfloor$ and $D=\varepsilon 2^m13^j17^k$ with $j, k\in\{0,1\}$.
As before, multiplying both sides of \eqref{eqp9} by $7^2D^3$, we get 
\begin{equation}\label{eqp10}
V^2=U^3-35dDU+147d^2D^2U-49d^3D^3,
\end{equation}
where $U=7DX$ and $V=7D^2Y$. Here, we use \texttt{SIntegralPoints} subroutine of MAGMA to compute all $\{5\}$-integral points on the $32$ elliptic curves defined by \eqref{eqp9}. As in the previous case, we apply the relations $U=7DX, V=7D^2Y, X=u^2/5^{2a_1}, Y=2^{m}13^{b_2}17^{c_2}/5^{3a_1}, a_2=\lfloor{a_1/2}\rfloor, b_2=\lfloor{b_1/2}\rfloor, c_2=\lfloor{c_1/2}\rfloor$ and $D=\varepsilon 2^m13^j17^k$ with $j, k\in\{0,1\}$, to find the corresponding integer solutions of \eqref{eqp1}. However,  none of these integral points leads to an integer solution of \eqref{eqp1}.

Assume that $v=17^{c_1}$ with $c_1\geq 1$. Then \eqref{eqp6} becomes 
$$2^{3m}5^{a_1}13^{b_1}=\varepsilon(7u^6-35u^4v^2d+21u^2v^4d^2-v^6d^3).
$$
As before, we first divide both sides of this equation by $v^6$ and then multiply by  $7^2D^3$ to arrive at
\begin{equation}\label{eqp11}
V^2=U^3-35dDU+147d^2D^2U-49d^3D^3,
\end{equation}
where $U=7DX=7Du^2/v^2, V=7D^22^{m}5^{a_2}13^{b_2}/v^3, a_2=\lfloor{a_1/2}\rfloor, b_2=\lfloor{b_1/2}\rfloor$ and $D=\varepsilon 2^m5^i13^j$ with $i,j\in\{0,1\}$. As in the previous case, we compute all $\{17\}$-integral points on the elliptic curves defined by \eqref{eqp10}, but  none of these integral points leads to an integral solution of \eqref{eqp1}.

Finally let $v=5^{a_1}17^{c_1}$ with $a_1,c_1\geq 1$. Then \eqref{eqp6} becomes 
$$2^{3m}13^{b_1}=\varepsilon(7u^6-35u^4v^2d+21u^2v^4d^2-v^6d^3).
$$
In the same way as before, we divide both sides of this equation by $v^6$ and then multiply by  $7^2D^3$ to get
\begin{equation}\label{eqp12}
V^2=U^3-35dDU+147d^2D^2U-49d^3D^3,
\end{equation}
where $U=7DX=7Du^2/v^2, V=7D^22^{m}13^{b_2}/v^3, b_2=\lfloor{b_1/2}\rfloor$ and $D=\varepsilon 2^m13^j$ with $j\in\{0,1\}$. We compute all $\{5,17\}$-integral points on the $16$ elliptic curves defined by \eqref{eqp11}, but  none of these integral points leads to an integer  solution of \eqref{eqp1}.

\subsection*{Case III: When $p=5$} 
Let $q$ be  a primitive divisor of $\mathfrak{L}_5(\alpha, \bar{\alpha})$. Then by \eqref{eqp4} only possibility $q=5$ or $q=13$ or $17$. Again the fact that any primitive divisor of $\mathfrak{L}_p(\alpha, \bar{\alpha})$ is congruent to $\pm 1$ modulo $p$ confirms that none of $5,13$ and $17$ is a primitive divisor of $\mathfrak{L}_5(\alpha, \bar{\alpha})$.  Therefore as in Case I, using Primitive Divisor Theorem for Lehmer sequences and  \cite[Table 4]{BH01}, we get 
\begin{equation}\label{eqp13}
(2^{2-m}u^2, 2^{2-m}v^2d)=\begin{cases}
(F_{k-2\varepsilon}, 4F_\varepsilon-F_{k-2\varepsilon}) \text{ with } k\geq 3,\\
(L_{k-2\varepsilon}, 4L_\varepsilon-L_{k-2\varepsilon}) \text{ with } k\ne 1,
\end{cases}
\end{equation}
where $F_k$ (resp. $L_k$) denotes the $k$-th Fibonacci (resp. Lucas) number. Utilizing Lemma \ref{lemCO} in \eqref{eqp11}, we conclude the following:
\begin{itemize}
\item $(k-2\varepsilon, m, 2u)=(1,0,1), (2,0,1), (12,0,12)$; none of these are possible as $u$ is odd.
\item $(k-2\varepsilon, m, u)=(3,1,1), (6,1,2)$, but the only possibility $(k-2\varepsilon, u)=(3,1,1)$. This leads to $(k, m,u,v,d)=(5,1,1,3,1)$, and thus $y=(u^2+dv^2)/2=5$. Therefore \eqref{eqp1} becomes $x^2+1=2\times 5^5$, which has no integer solutions.  
\item $(k-2\varepsilon, m, 2u)=(1,0,1), (3,0,2)$, but the only possibility is $(k-2\varepsilon, m, u)=(3,0,1)$. This leads to $v^2d=8$, which is not possible as $vd$ is odd.
\item $(k-2\varepsilon, m, u)=(6,1,3)$, which gives $(k,m, u, vd)=(4,1,3,5), (8,1,3,89)$. The only possibility is $(k,m, u, vd)=(4,1,3,5)$, which yields $y=7$. Thus \eqref{eqp1} becomes $x^2+5=2\times 7^5$, which has no integer solutions.  
\end{itemize}

\subsection*{Case IV: $p=3$}  
In this case, the facts of primitive divisors of $\mathfrak{L}_3(\alpha, \bar{\alpha})$ does not provide any fruitful conclusion. Thus, we transfer the problem of finding integer solutions of \eqref{eqp1} into the problem of finding $\{5,13,17\}$-integral points on the corresponding elliptic curves. For $p=3$, \eqref{eqp1} becomes
\begin{equation}\label{eqp14}
x^2+5^a13^b17^c=2^my^3,~ x,y\geq 1, \gcd(x, y)=1, a,b,c\geq 0, p\geq 2, m=0,1.
\end{equation} 
We write 
\begin{equation}\label{eqp15}
2^{2m}5^a13^b17^c=Az^6,
\end{equation}
where $(A,z)=\left(2^{2m}5^{a_1}13^{b_1}17^{c_1}, 5^{a_2}3^{b_2}7^{c_2}\right)$ with $a=6a_2+a_1, b=6b_2+b_1, c=6c_2+c_1$,  $a_1,b_1,c_1\in\{0,1,2,3,4,5\}$ and $a_2,b_2, c_2\geq 0$.

We now multiply the both sides of \eqref{eqp14} by $2^{2m}$ and then subsequently divide both sides of $z^6$, and put $X:=2^mx/z^3$ and $2^my/z^2$ to reduce to 
\begin{equation}\label{eqp16}
X^2=Y^3-A.
\end{equation}
Here, we again use \texttt{SIntegralPoints} subroutine of MAGMA  to determine all $ \{5,13,17\}$-integral points on the $432$ elliptic curves defined by \eqref{eqp16}. Taking into account that $xy \ne 0$ and $\gcd(x, y) = 1$, we don't consider the $\{5,13,17\}$-integer points such that $XY=0$ or the numerators of $X/2^m$ and $Y/2^m$ are not coprime. Finally, we utilize the relation \eqref{eqp15} along with the conditions on $A$ and $z$, and on their exponents to find the corresponding integer solutions of \eqref{eqp14}. These solutions are given by  $(x,y,a,b,c,m)\in\mathfrak{S}$.
\end{proof}

\section{Proof of Theorem \ref{thm}}\label{S4}
We first assume that  $n$ is a multiple of $4$. Then all integer solution of \eqref{eqn1} are given by Proposition \ref{propm4} (see, Table \ref{Tm4}). 

Let $n\geq 3$ be an integer such that $4\nmid n$. Then we can write $n=p\ell$ for some odd prime $p$ and some odd integer $\ell\geq 1$. Thus, \eqref{eqn1} can rewritten as 
\begin{equation}\label{eqnf1}
x^2+5^a13^b17^c=2^m (y^{n/p})^p,~~ x\geq 1, y>1, \gcd(x, y)=1, a, b,c, m\geq 0, p\geq 3.
\end{equation}   
By Proposition \ref{propp}, the only integer solutions of \eqref{eqnf1} are given by 
\begin{align*}
(x,y^{n/p},a,b,c,m)\in&\mathfrak{S}. 
\end{align*}

These solutions lead to the integer solutions of \eqref{eqn1}, which are listed in Table \ref{Tn}.
			
\begin{center}
{\small
\begin{longtable}{ccccccc| ccccccc}
\caption{All the solutions of \eqref{eqn1} when $3\mid n$} \label{Tn} \\

\hline \multicolumn{1}{c}{$x$} & \multicolumn{1}{c}{$y$} & \multicolumn{1}{c}{$a$}& \multicolumn{1}{c}{$b$} & \multicolumn{1}{c}{$c$}& \multicolumn{1}{c}{$m$} &\multicolumn{1}{c}{$n$} &\multicolumn{1}{|c}{$x$} & \multicolumn{1}{c}{$y$} &\multicolumn{1}{c}{$a$} &\multicolumn{1}{c}{$b$} & \multicolumn{1}{c}{$c$}& \multicolumn{1}{c}{$m$} & \multicolumn{1}{c}{$n$}\\ \hline 
\endfirsthead

\multicolumn{14}{c}%
{{\bfseries \tablename\ \thetable{} -- continued from previous page}} \\
\hline \multicolumn{1}{c}{$x$} & \multicolumn{1}{c}{$y$} & \multicolumn{1}{c}{$a$}& \multicolumn{1}{c}{$b$} & \multicolumn{1}{c}{$c$}& \multicolumn{1}{c}{$m$} &\multicolumn{1}{c}{$n$} &\multicolumn{1}{|c}{$x$} & \multicolumn{1}{c}{$y$} &\multicolumn{1}{c}{$a$} &\multicolumn{1}{c}{$b$} & \multicolumn{1}{c}{$c$}& \multicolumn{1}{c}{$m$} & \multicolumn{1}{c}{$n$}\\ \hline 
\endhead
\hline \multicolumn{14}{c}{{Continued on next page}} \\ \hline
\endfoot
\hline 
\endlastfoot
70&17&0&1&0&0&3&
 716&81&1&1&2&0&3\\
 716&9&1&1&2&0&6&
 716&3&1&1&2&0&12\\
 94&21&2&0&1&0&3&
 142&29&2&2&0&0&3\\
2034&161&3&0&2&0&3&
9&5&0&2&0&1&3\\
7&3&1&0&0&1&3&
99&17&2&0&0&1&3\\
63&13&2&0&1&1&3&
19&7&2&1&0&1&3 \\
33&7&2&2&1&1&3& 
118699& 1917&2&2&1&1&3\\
79137& 1463& 2&3&0&1& 3&
253&73&2&4&0&1&3 \\
188000497& 260473& 8&4&0&1&3&
267689& 3297& 2&2&3&1&3\\
336049& 4317& 10&0&3&1&3& 17127& 553& 6& 2& 1&1& 3
\end{longtable}}
\end{center}
\vspace{-15mm}
\section*{acknowledgements}
The authors are grateful to Professor Kalyan Chakraborty for carefully reading and valuable comments which have improved the presentation of this manuscript. The authors are thankful to the anonymous referee(s) for valuable comments/suggestions.  
This work is supported by SERB MATRICS Project (No. MTR/2021/000762) and SERB CRG Grant (No. CRG/2023/007323), Govt. of India.


\begin{thebibliography}{25}

\bibitem{AA99} F. S. Abu Muriefah and S. A. Arif, {\it The Diophantine equation $x^2 + 5^{2k+1} = y^n$}, Indian J. Pure Appl. Math. {\bf 30} (1999), no. 3, 229--231.

\bibitem{AL08} F. S. Abu Muriefah, F. Luca and A. Togb\'e, {\it  On the Diophantine equation $x^2 + 5^a 13^b = y^n$},  Glasgow Math. J. {\bf 50} (2008), no. 1, 175--181.
\bibitem{AL09} F. S. Abu Muriefah, F. Luca,  S. Siksek and S. Tengely, {\it On the Diophantine equation $x^2 + C = 2y^n$}, {\it Int. J. Number Theory} {\bf 5} (2009), no. 6,  1117--1128.

\bibitem{AZ20} M. Alan and U. Zengin, {\it On the Diophantine equation $x^2 + 3^a41^b = y^n$}, Period. Math. Hung. {\bf 81} (2020), 284--291.

\bibitem{AA02} S. A. Arif and A. S. Al-Ali, {\it On the Diophantine equation $ax^2+b^m=4y^n$},  Acta Arith. {\bf 103} (2002), no. 4, 343--346.

\bibitem{BHS19} S. Bhatter, A. Hoque and R. Sharma, {\it On the solutions of a Lebesgue-Nagell type equation},  Acta Math. Hungar. {\bf 158} (2019), 17--26.

\bibitem{BH01} Y. Bilu, G. Hanrot and P. M. Voutier, {\it Existence of primitive divisors of Lucas and Lehmer numbers (with an appendix by M. Mignotte)}, J. Reine Angew. Math. {\bf 539} (2001), 75--122.

\bibitem{BCP97} W. Bosma, J. Cannon and C. Playoust, {\it The Magma algebra system. I. The user language}, J. Symbolic Comput. {\bf 24} (1997), no. 3, 235--265.

\bibitem{BU01} Y. Bugeaud, {\it On some exponential Diophantine equations}, Monatsh. Math. {\bf 132} (2001), 93--97.


\bibitem{CH22} K. Chakraborty and A. Hoque, {\it On the Diophantine equation $dx^2+p^{2a}q^{2b}=4y^p$}, Results Math. {\bf 77} (2022), no. 1, pp 11, article no. 18.


\bibitem{CHS20} K. Chakraborty, A. Hoque and R. Sharma, {\it Complete solutions of certain Lebesgue-Ramanujan-Nagell equations}, Publ. Math. Debrecen, {\bf 97} (2020), no. 3/4, 339--352.

\bibitem{CHS-21} K. Chakraborty, A. Hoque and R. Sharma, {\it On the solutions of certain Lebesgue-Ramanujan-Nagell equations}, Rocky Mountain J. Math. {\bf 51} (2021), no. 2, 459--471.

\bibitem{CHS21} K. Chakraborty, A. Hoque and K. Srinivas, {\it On the Diophantine equation $cx^2+p^{2m}=4y^n$}, Results Math. {\bf 76} (2021),  pp. 12, article no. 57.

\bibitem{CO64}  J. H. E. Cohn, {\it Square Fibonacci numbers, etc.}, Fibonacci Quart. {\bf 2} (1964), no. 2, 109--113.

\bibitem{DA11} A. Dabrowski, {\it On the Lebesgue-Nagell equation}, Colloq. Math. {\bf 125} (2011), no. 2, 245--253.

\bibitem{DGS21} A. Dabrowski, N. G\"unhan and G. Soydan, {\it On a class of Lebesgue-Ljunggren-Nagell type equations}, J. Number Theory {\bf 215} (2020), 149--159.

\bibitem{DE17} M. Demirci, {\it On the Diophantine equation $x^2 + 5^a p^b = y^n$}, Filomat {\bf 31} (2017), no. 16, 5263--5269.


\bibitem{GW12} S. Gou and T. T. Wang, {\it The Diophantine equation $x^2 + 2^a.17^b = y^n$}, Czechoslovak Math. J. {\bf 62} (2012), 645--654.

\bibitem{HO20} A. Hoque, {\it On a class of Lebesgue-Ramanujan-Nagell equations}, Period. Math. Hungar.   \url{https://doi.org/10.1007/s10998-023-00564-z}



\bibitem{LS20} M. Le and G. Soydan, {\it A brief survey on the generalized Lebesgue-Ramanujan-Nagell equation},  Surveys in Mathematics and its Applications {\bf 15} (2020), 473--523.

\bibitem{LE50} V. A. Lebesgue, {\it Sur l'impossibilit\'e, en nombres entiers, de l'\'equation $x^m =y^2+1$}, Nouvelles Annales des Math. {\bf 9} (1850), p. 178.

\bibitem{LJ66} W. Ljunggren, {\it  On the Diophantine equation $Cx^2 + D = 2y^n$},  Math. Scand. {\bf 18} (1966), 69--86.


\bibitem{LT09} F. Luca and A. Togb\'e, {\it On the equation $x^2 + 2^\alpha 13^\beta = y^n$}, {\it Colloq. Math.} {\bf 116} (2009), no. 1, 139--146.

\bibitem{PR11}  I. Pink and Z. R\'abai, {\it On the Diophantine equation $x^2 + 5^k17^l = y^n$}, Commun. Math. {\bf 19} (2011), 1--9. 



\bibitem{TA09} L. Tao, {\it On the Diophantine equation $x^2 + 5^m = y^n$}, Ramanujan J. {\bf 19} (2009), 325--338.

\bibitem{YU05} P.  Yuan, {\it On the Diophantine equation $ax^2 + by^2 = ck^n$}, Indag. Math., N. S.  {\bf 16} (2005), no. 2, 301--320.

\bibitem{ZL15} H. Zhu, M. Le, G. Soydan and A. Togb\'e, {\it On the exponential Diophantine equation $x +2^a p^b = y^n$},  Period. Math. Hung. {\bf 70} (2015), 233--247.

\end{thebibliography}
\end{document}